\documentclass[a4paper]{article}

\usepackage{graphics}

\usepackage{graphicx}
\usepackage{caption}
\usepackage{subcaption}

 \usepackage{algorithm}
 \usepackage{algpseudocode}

\algnewcommand{\Multcomment}[1]{\\\hfill \begin{minipage}[t]{0.8\textwidth}/* {\itshape #1} */\end{minipage}}

\algnewcommand{\algorithmicgoto}{\textbf{go to}}%
\algnewcommand{\Goto}[1]{\algorithmicgoto~\ref{#1}}%

\usepackage{enumitem} 
\usepackage{amsmath}
\usepackage{amssymb}
\usepackage{graphicx}
\usepackage[colorinlistoftodos]{todonotes}
\usepackage[colorlinks=true, allcolors=blue]{hyperref}
\usepackage{color}
\usepackage{multirow}
\usepackage[round, sort,comma,authoryear]{natbib}
\usepackage{siunitx} 
\usepackage{verbatim}
\usepackage{breakurl}

\usepackage{tikz}

\hypersetup{colorlinks, citecolor=black, linkcolor=black, urlcolor=black}

\definecolor{brightturquoise}{rgb}{0.03, 0.91, 0.87}

\definecolor{mygreen}{RGB}{0,128,0}

\def\rev#1{{\color{black}{#1}}} 

\newcommand{\footnoteremember}[2]{
\footnote{#2}
  \newcounter{#1}
  \setcounter{#1}{\value{footnote}}
}
\newcommand{\footnoterecall}[1]{
\footnotemark[\value{#1}]
}

\title{Pseudo-random Instance Generators in C++ for Deterministic and Stochastic Multi-commodity Network Design Problems}

\author{
Eric Larsen\footnoteremember{cirrelt-UdeM}{CIRRELT and Universit\'e de Montr\'eal}
\and
Serge Bisaillon\footnoterecall{cirrelt-UdeM}
\and
Jean-Fran\c{c}ois Cordeau\footnoteremember{cirrelt-HEC}{CIRRELT and HEC Montr\'eal}
\and
Emma Frejinger\footnoterecall{cirrelt-UdeM}\footnote{Corresponding author. Email: emma.frejinger@umontreal.ca}
}

\date{\today}

\begin{document}

\maketitle

\begin{abstract}Network design problems constitute an important family of combinatorial optimization problems for which many exact and heuristic algorithms have been developed over the last few decades. Two central problems in this family are the multi-commodity, capacitated, fixed charge network design problem (MCFNDP) and its stochastic counterpart, the two-stage stochastic MCFNDP with recourse. These are standard problems that often serve as workbenches for devising and testing models and algorithms in stylized but near-realistic settings. Progress on solution methods requires instances of varying size, structure and difficulty. We present two open-source C++ generators. The deterministic generator modernizes and extends the generator underlying the Canad instances, with new topologies and continuous control of capacity tightness and cost structure. The stochastic generator, which is to our knowledge the first publicly available one for the two-stage stochastic MCFNDP, builds scenarios that match user-specified moments and correlations. We also release libraries of configurations for 6,000 deterministic and 416 stochastic instances, as well as the base instances and configurations of 20 stochastic instances on hub-and-spoke benchmark networks. The generators aim to support systematic and reproducible experimentation, including the training and out-of-sample evaluation of learning-based and stochastic programming methods.
\end{abstract}


\section{Introduction}

Network design problems \citep{CraiEtAl2021a} are concerned with choosing which arcs to open in a  network and how to route flows from origins to destinations so as to satisfy a set of demands while respecting arc capacities. These flows can represent objects, people or information, and may be finely or coarsely differentiated depending on the application. Most variants of network design can be characterized by the following four attributes: (i)~the amounts of commodities demanded and supplied at the nodes of the network, (ii) the fixed costs associated with the opening of arcs, (iii) the flow-dependent and commodity-specific costs on the arcs, and (iv) the global and commodity-specific capacities on the arcs.

NP-hard network design problems are ubiquitous and occur in many areas of human activity such as supply chain management, transportation, and telecommunications. \rev{Developing exact and heuristic solution methods for these problems is a very active field of research in combinatorial optimization.} In this context, the linear, deterministic, multi-commodity, capacitated, fixed charge network design problem (MCFNDP) \citep{gendron1999multicommodity} and its stochastic counterpart, the two-stage stochastic MCFNDP with recourse \citep{crainic2011progressive}, play a central role. They often serve as \rev{workbenches} for the design and testing of new mathematical formulations and solution algorithms. In addition, the successes achieved on these two fundamental problems often carry over to other related problems such as facility location problems. \rev{Progress on these methods requires instances of varying difficulty and structure for testing and benchmarking. Our objective is a software contribution providing the community with the means to produce such instances systematically and reproducibly. Our focus is on the core features of the standard MCFNDP, which is a cornerstone of the network design literature, and not on specific applications. However, the two generators are flexible enough to be tailored to most applications.}

\rev{Benchmark instances are ubiquitous in the operations research community for comparing algorithm performance, a prominent example being the MIPLIB library \citep{GleixnerEtAl21}. For the deterministic MCFNDP, the Canad instances \citep{Canad} have served this role for over three decades, and they have been used to generate limited sets of stochastic instances \citep[e.g.,][]{RahmanianiEtAl24}. While such fixed instance sets are invaluable for comparing algorithms, relying on them for algorithm \emph{development} carries a risk of over-tuning to their idiosyncrasies, especially if the set is small or lacks variability in instance sizes, difficulty, or structure. This concern is compounded in two settings motivating our contribution. \emph{First}, the growing use of machine-learning-informed algorithms requires large numbers of instances, distinct from the test set, for training and validation \citep[][discuss related issues and introduce ``Distributional MIPLIB'' for testing]{huangEtAl2024}. \emph{Second}, evaluation of stochastic programming methods calls for generating many instances, both to assess the stability of a scenario-generation method and to evaluate solutions out of sample \citep[see][for a survey on sampling-based methods for stochastic optimization]{HomeBayr2014}. Both settings call for access to a distribution of instances from which many can be drawn. Instance generators serve this purpose.} To the best of our knowledge, there does not exist a publicly available generator for the two-stage stochastic MCFNDP. 

\rev{A published and documented generator also supports reproducibility, and an adaptable one makes it possible to build test sets that are deliberately diverse, to support general conclusions, or specialized, to support conclusions about a given class of instances. This paper makes the following infrastructure contributions:
\begin{itemize}
\item \rev{We introduce what is, to our knowledge, the first open-source instance generator for the two-stage stochastic MCFNDP with recourse. Published studies have so far randomized deterministic instances ad hoc, often with the moment-matching algorithm of \citet{HoylEtAl2003}, hereafter HKW, and instance-specific targets \citep[see, e.g.,][]{crainic2011progressive,CraiEtAl2021b}. Our generator derives the target moments from a uniform or triangular distribution around each base value and the target correlations from a few block values, so that a whole experiment is specified by a handful of parameters and a pair of pseudo-random seed and stream. Any base instance can be used, and every scenario can be checked for feasibility and validity. The generator thereby supports the systematic and reproducible generation of stochastic instances.}
\item The stochastic generator takes a deterministic MCFNDP instance as input, which supplies the network topology and a base scenario. This instance can be provided by the user in the required input format, or generated pseudo-randomly with our modernized and extended version of Mulgen \citep{Mulgen}, the generator underlying the Canad instances \citep{Canad} that are widely used \citep[see, e.g.,][]{Chouman2016,Crainic2001,Hewitt2009}. The revised generator also stands on its own as a generator of deterministic MCFNDPs.
\item Alongside the generators, we release two libraries of configurations from which benchmark instances can be regenerated exactly: (i) 6,000 configurations of deterministic instances, spanning two primary topologies, four secondary topologies, two origin–destination structures, networks of 6 to 476 nodes and 3 to 312 commodities, and a range of capacity tightness and cost balance; (ii) 416 configurations of stochastic instances built on 32 base instances. These libraries are a starting point rather than a fixed benchmark: as discussed in Section~\ref{sec:det_oper}, any fixed set eventually invites over-tuning, and the generators allow the libraries to be extended as methods and experimental needs evolve.
\end{itemize}}

\rev{Pseudo-randomly generated instances are valuable for producing instance distributions with varying characteristics. Users seeking real-world or realistic network topologies for a specific application domain, such as telecommunications, transportation or supply chain management, can instead supply a user-defined topology. Other generators may also serve this complementary purpose. For instance, the generator of \citet{BonnEtAl2025} produces deterministic MCFNDPs and schedule-augmented MCFNDPs (SMCFNDPs)} following a top-down approach: the user sets targets for metrics taken from the literature on complex network analysis, and the generator produces instances with a hub-and-spoke topology in a process where the user aims to closely match these targets. Our generator instead follows a bottom-up approach: the user sets elementary characteristics of the MCFNDP, together with the values or ranges its individual elements may take, and the generator selects these elements pseudo-randomly within the specified ranges.

\rev{The remainder of the paper is organized as follows. Section~\ref{sec:det_gen} recalls the deterministic MCFNDP and describes the deterministic generator, its workflow, its library and its performance. Section~\ref{sec:stoch_gen} introduces the two-stage stochastic MCFNDP with recourse, describes the stochastic generator and its library, and illustrates its use on the networks of \citet{BonnEtAl2025}. Section~\ref{sec:concl} concludes.}

\section{A generator of linear deterministic MCFNDPs} \label{sec:det_gen}
This section first introduces the linear deterministic MCFNDP \citep[borrowing notation from][]{CraiEtAl2021a} and follows with a description of the corresponding deterministic generator.

\subsection{The linear deterministic MCFNDP} \label{sec:det_prob}

Let $\mathcal{G}=(\mathcal{N},\mathcal{A})$ denote a graph composed of arcs $(i,j)\in \mathcal{A}$ and nodes $i,j\in\mathcal{N}$. \rev{A fixed cost $f_{ij}$ and a capacity $u_{ij} \ge 0$ limiting the total amount of flow on the arc are associated with each arc $(i,j) \in \mathcal{A}$}. Demands are associated with commodities $k \in \mathcal{K}$ and defined over the nodes of the graph. The unit costs for using any arc $(i,j) \in \mathcal{A}$ are denoted $c_{ij}^k$. It is commonly assumed that each commodity is characterized by a single origin node $O(k)\in \mathcal{N}$, a single destination node $D(k)\in\mathcal{N}$ and a volume $d^k$.  In this setting, the net outgoing flow of commodity $k$ at node $i$ is given by
\[
w_i^k =
\begin{cases}
d^k, & \text{if } i = O(k), \\
-d^k, & \text{if } i = D(k), \\
0, & \text{otherwise}.
\end{cases}
\]
The problem may also include commodity-specific capacities $b_{ij}^k\geq 0$ limiting the flow of specific commodities on the arcs. Note \rev{however} that the generator introduced below allows the \rev{existence} of multiple sources and sinks for each commodity.

The linear deterministic MCFNDP comprises two sets of decision variables -- binary design variables \rev{$y_{ij}\in \{0,1\},~(i,j)\in\mathcal{A}$}, and continuous multi-commodity flow variables $x_{ij}^k\geq 0,~(i,j)\in\mathcal{A},k\in\mathcal{K}$. The problem can be defined as
\begin{equation}
\label{eq:glob_obj_det} \min_{\boldsymbol{y},\boldsymbol{x}}~\bigg\lbrace \sum_{(i,j) \in \mathcal{A}} f_{ij} y_{ij} + \sum_{k\in \mathcal{K}} \sum_{(i,j) \in \mathcal{A}} c^k_{ij} x^k_{ij} \bigg\rbrace \\ 
\end{equation}

\begin{alignat}{3}
\label{eq:glob_const1_det} \text{subject to } & \sum_{j \in \mathcal{N}_i^+} x_{ij}^k - \sum_{j \in \mathcal{N}_i^-} x_{ji}^k = w_i^k, & \forall i \in \mathcal{N}, \forall k \in \mathcal{K}, \\
\label{eq:glob_const2_det}  & \sum_{k \in \mathcal{K}} x_{ij}^k \leq u_{ij} y_{ij}, & \forall (i,j) \in \mathcal{A}, \\
\label{eq:glob_const2f_det}  & x_{ij}^k \leq b_{ij}^k y_{ij}, & \forall (i,j) \in \mathcal{A}, \forall k \in \mathcal{K}, \\
\label{eq:glob_const3_det} & x_{ij}^k \geq 0, & \forall (i,j) \in \mathcal{A}, \forall k \in \mathcal{K}, \\
\label{eq:glob_const4_det} & y_{ij} \in \lbrace 0,1 \rbrace, & \forall (i,j) \in \mathcal{A},
\end{alignat} 
where $\boldsymbol{y}$ \rev{and} $\boldsymbol{x}$ respectively stand for all design variables $y_{ij},~(i,j)\in\mathcal{A}$, and all flow variables $x_{ij}^k,~(i,j)\in\mathcal{A},k\in\mathcal{K}$. The objective function (\ref{eq:glob_obj_det}) minimizes the total costs expressed as a summation of fixed arc opening costs and a summation of flow-dependent transportation costs. Constraints~(\ref{eq:glob_const1_det}) enforce flow conservation at each node $i\in \mathcal{N}$, where $\mathcal{N}_i^+$ and $\mathcal{N}_i^-$ identify, respectively, the successor and predecessor nodes of $i$.
Constraints~(\ref{eq:glob_const2_det}) enforce overall capacity limits and act as linking constraints. Finally, \rev{Constraints}~(\ref{eq:glob_const2f_det}) \rev{are} commodity-specific capacity constraints (a.k.a. strong forcing constraints) restricting the flow of commodity $k\in\mathcal{K}$ on arc $(i,j)$. The latter are optional in the generator.

\subsection{The generator} \label{sec:det_gen_tool}

\rev{We first review the distinctive characteristics of the generator (Section~\ref{sec:det_design}), then outline its workflow (Section~\ref{sec:det_workflow}), and finally describe its operation, the accompanying library and its performance (Section~\ref{sec:det_oper}).}

\subsubsection{Design.} \label{sec:det_design}

\rev{The revised generator differs from its predecessor in four respects. \emph{First}, topologies are more varied: (i) the graph can be read from or written to a text file, so that bespoke or real topologies can be used, (ii) hubs-and-satellites primary topologies (of which hub-and-spoke is a particular case) complement the point-to-point topology of the original generator, (iii) grid-like topologies can lie on cylinders (as in cyclical schedules in space-time settings) or on flat surfaces (as in bounded networks), and not only on a torus. \emph{Second}, the tightness between overall volumes and overall capacities, and the balance between overall fixed and variable costs, can be adjusted continuously in both directions. \emph{Third}, the generator guides the user through automatic safeguards, can draw small instances in TeX format, and is extensively documented. \emph{Fourth}, an auxiliary program verifies the feasibility of the generated instances. The C++ code also replaces the original pseudo-random number generator with the PCG library \citep[see][]{PCGSite, LemireSite}, which is faster, statistically more reliable and supports random streams, so that each instance can be retrieved from its seed and stream.}

\subsubsection{\rev{Workflow.}} \label{sec:det_workflow}

\rev{Generation proceeds through the following abstract steps.}
\begin{enumerate}
\item \rev{Determine the primary topology: either point-to-point or hubs-and-satellites.}
\item \rev{Determine the secondary topology: either grid-like (toral, cylindrical or flat), circular or random. Under point-to-point, the secondary topology determines all nodes and arcs. Under hubs-and-satellites, the secondary topology characterizes only hub nodes and arcs between hubs whereas each hub is assigned a set of satellite nodes connected sequentially and arcs connect the hub to both ends of the sequence.}
\item \rev{Generate commodity volumes in one of three ways: (i) single source node and single sink node for each commodity, (ii) multiple common source nodes for all commodities and multiple common sink nodes for all commodities, (iii) random source nodes and sink nodes for each commodity. Under point-to-point, any node can be a source or a sink of volume. Under hubs-and-satellites, only satellite nodes can be sources or sinks of volumes. When the secondary topology is grid-like or circular, Steps 1 and 2 ensure that the graph is strongly connected and that the non-random arcs have capacities sufficient to carry all commodity volumes. No such guarantee is given for a random secondary topology, nor for a graph read from a file, whose connectivity is the user’s responsibility.}
\item \rev{Introduce additional, random arcs with lesser capacities and lesser costs, thus defining new and more advantageous paths from sources to sinks.}
\item \rev{Adjust the commodity-specific capacities and fixed costs of the random arcs of Step 4, which makes these shortcuts more or less advantageous and thus changes the difficulty of the design problem.}
\item \rev{Adjust the overall tightness between volumes and capacities and the overall balance between fixed and variable costs, upward or downward, around the neutral value 1. These two adjustments generalize the classical distinction, used to build the Canad instances, between instances dominated by fixed or by variable costs and between instances with loose or tight capacities, and they are the main levers for diversifying a test set. Raising fixed costs relative to variable costs makes design decisions dominate the objective and typically weakens the linear relaxation, whereas lowering them brings the problem closer to a multi-commodity flow problem. Tightening capacities increases the interaction between commodities on the arcs and tends to increase difficulty, until instances become infeasible. The auxiliary program detects infeasible instances. Instances dominated by fixed costs and with tight capacities are usually the most challenging. Marginally feasible instances are not suitable for randomization, since some of their scenarios would be infeasible: a range of tightness settings should be used when instances are meant to be randomized. Implementing ranges of values for the adjustments of capacity tightness and overall balance between costs also serves the goal of diversification. For example, in the library being released a typical range for both adjustments is [0.33, 0.5, 1, 2, 3].}
\item \rev{Write the instance in the two MCFNDP-specific formats DOW and STD, and in the generic LP and MPS formats and, optionally, draw it in TeX format. The drawing shows the graph, its legend and, when costs or capacities differ between commodities, a third pane (Figure 1 shows the main pane for a 3 x 4 flat grid with three commodities). It helps the user explore the options of the generator and choose the characteristics of the networks. Since such drawings quickly become unreadable and costly to compile, the generator produces them only for small instances. The user guide states the limits.}
\end{enumerate}

\subsubsection{\rev{Operation, library} and performance.} \label{sec:det_oper}

\rev{The generator, its user guide (USERGUIDE\_DETERM.md) and an extensively commented tutorial example (TUTOREXAMPLE\_DETERM.md) are available at \url{https://github.com/larseeri/determ-MCFNDP-gener}. The generator can be built from source and run on Linux, or run through the supplied Docker image.}

\rev{The Canad instances, created more than 30 years ago with the original generator, are still a main test set for the MCFNDP. They have three limitations. \emph{First}, the number of Canad instances is small (e.g., there are 162 Canad R instances), which means that the same instances must be reused continually in testing experiments. Continual reuse is highly vulnerable to over-tuning and hence to a progressive loss of evidential power. \emph{Second}, computational methods have progressed immensely over the last 30 years and MCFNDP instances that were reasonably challenging at their introduction have ceased to be. For instance, the Gurobi solver \rev{version 12.0.0 running on an AMD EPYC 7763 processor} can respectively solve to optimality 60\%, 70\%, 80\% and 90\% of all Canad R instances within less than 10 s, 56 s, 744 s and 21750 s with a single core and default settings. For only 13 instances the limit of 10 h is exceeded and even in those cases the final optimality gap is less than 3.75\%. \emph{Third}, Canad instances lack structural diversity: they all present a point-to-point grid topology on a torus, commodities with a single source and a single sink, variable costs that do not depend on the commodity, and commodity-specific capacities that are simply equated to the minimum of overall capacity and commodity volume. A generator addresses these limitations: it can produce a diversified set of instances and extend it over time.}

\rev{We created a library of 6,000 configurations of deterministic MCFNDP
instances. To save space, only the configuration files are stored in the /library folder of the online repository. Each instance is regenerated in seconds from its configuration, which fixes the seed and stream of the pseudo-random generator. The library comprises 3,600 point-to-point and 2,400 hubs-and-satellites configurations, with toral, cylindrical, flat or circular secondary topologies, single or common sources and sinks, networks of 6 to 476 nodes (6 to 72 hub nodes) and 3 to 312 commodities, and a grid of capacity-tightness and cost-balance adjustments around the neutral value 1. The document DESCLIBRARY\_DETERM.txt available at this location details its content. The library is meant as a broadly diversified starting set, from which users can select or sample the instances relevant to their experiments, and which can be extended as methods and needs evolve.}

The generator achieves high speed computations. For example, generating and saving instances similar to the largest ones among the historical R and C series available at \citet{Canad}, that is r18 (20 nodes, 315 arcs, 200 commodities) and c64 (30 nodes, 700 arcs, 400 commodities), requires respectively less than 1 second and less than 4 seconds, using a single core of an Intel(R) Xeon(R) E5-2637 v4 @ 3.50GHz CPU running on AlmaLinux 9.2. The advanced user interested in the inner workings of the code will find \rev{it} useful that it has been extensively commented and documented.

\begin{figure}[htbp] 
    \centering
    \includegraphics[width=0.6\textwidth]{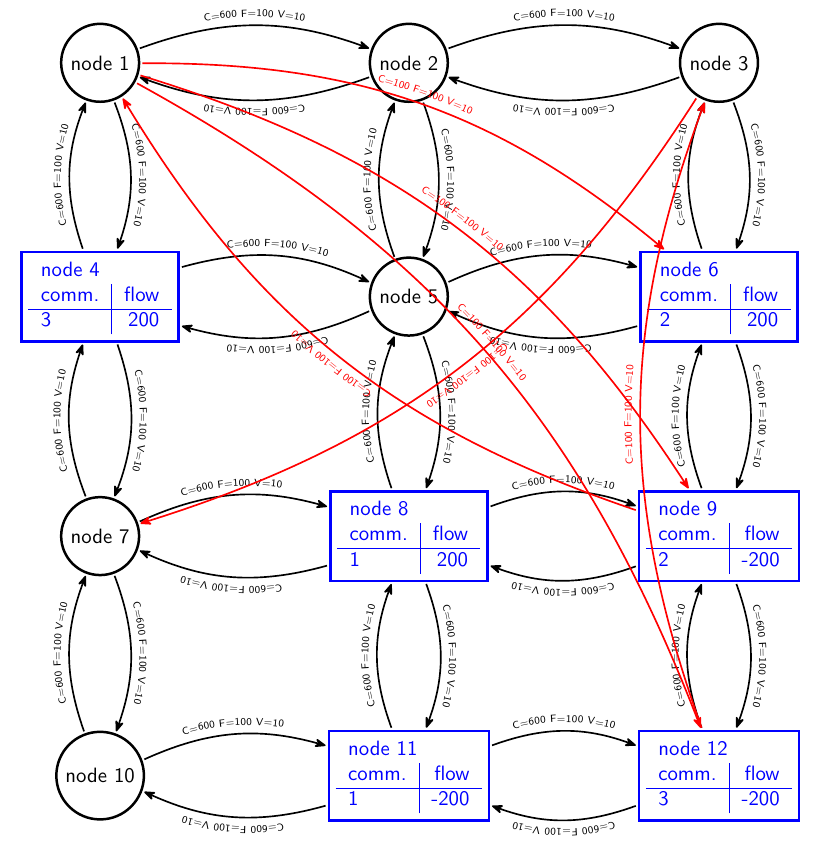} 
    \caption{Main pane of a small-scale example of graphical representation\\
    \tiny Reproduced from the legend pane:\\ Black arcs are part of the base, grid topology; red arcs are additional, randomly generated; $C=x$ along arc: overall capacity is equal to x; $F=y$ along arc: fixed cost of arc equal to y; $V=z$ along arc: variable cost of arc equal to z; capacities and variable costs are not differentiated between commodities; circular nodes are transit nodes; rectangular nodes are source ($\text{flow}>0$) or sink ($\text{flow}<0$) nodes for stated commodities.}
  
    \label{fig:example}
\end{figure}

\section[A generator of linear two-stage stochastic MCFNDPs with recourse]{A generator of linear two-stage stochastic \\ MCFNDPs with recourse} \label{sec:stoch_gen}

\rev{The stochastic generator randomizes a user-specified set of parameters of a supplied instance of the linear deterministic MCFNDP. It synthesizes their joint probability distribution, which is discrete for tractability: a finite set of realizations (scenarios) with their probabilities.}

\rev{This distribution is shaped by simple, interpretable requirements set by the user according to the questions under study, and by pseudo-random draws. The base instance can be real, hand-built or synthetic, produced for instance by the generator of Section~\ref{sec:det_gen} or by any other generator.}

\rev{Scenario-construction methods from the stochastic programming literature (sampling, reduction, generation) approximate a distribution that is known, or that can be sampled or estimated from data, so as to make the problem tractable. Our task differs: we synthesize a distribution that is not known a priori, from user requirements and pseudo-random draws. Appendix A reviews these methods and explains why we build scenarios with the Høyland-Kaut-Wallace (HKW) algorithm \citep{HoylEtAl2003}, which matches targets for the first four moments and the correlations of all random parameters. Section~\ref{sec:stoch_design_workflow} shows how these targets are derived from a few assumptions.}

\subsection{The linear two-stage stochastic MCFNDP with recourse} \label{sec:sto_prob}

\rev{We now introduce the notation for t}he linear two-stage stochastic MCFNDP with recourse. In comparison with the linear deterministic MCFNDP of Section~\ref{sec:det_prob}, $\omega$, $p(\omega)$ and $\Omega$ respectively identify a stochastic realization (i.e., a scenario), the probability of a realization and the set of all stochastic realizations. The latter is assumed to be finite. 
We also assume that parameters $\boldsymbol{w}$, $\boldsymbol{u}$, $\boldsymbol{b}$, and $\boldsymbol{c}$ may all have been randomized. Decision variables $\boldsymbol{y}$ and $\boldsymbol{x}(\omega)$ are respectively selected before and after disclosure of the stochastic realization $\omega$. The formulation becomes the following, where we have assumed for simplicity as in Section~\ref{sec:det_prob} that there is a single source node and a single sink node for each commodity:

\begin{equation}
\label{eq:glob_obj_sto} \min_{\boldsymbol{y}, \boldsymbol{x}(\omega), \forall \omega \in \Omega}~\bigg\lbrace \sum_{(i,j) \in \mathcal{A}} f_{ij} y_{ij} + \sum_{\omega \in \Omega} p(\omega)  \sum_{k\in \mathcal{K}} \sum_{(i,j) \in \mathcal{A}} c^k_{ij}(\omega) x^k_{ij}(\omega) \bigg\rbrace \\ 
\end{equation}

\begin{alignat}{3}
\label{eq:glob_const1_sto} \text{subject to } & \sum_{j \in \mathcal{N}_i^+} x_{ij}^k(\omega) - \sum_{j \in \mathcal{N}_i^-} x_{ji}^k(\omega) = w_i^k(\omega), & \forall i \in \mathcal{N}, \forall k \in \mathcal{K}, \forall \omega \in \Omega, \\
\label{eq:glob_const2_sto}  & \sum_{k \in \mathcal{K}} x_{ij}^k(\omega) \leq u_{ij}(\omega) y_{ij}, & \forall (i,j) \in \mathcal{A}, \forall \omega \in \Omega, \\
\label{eq:glob_const2f_sto}  & x_{ij}^k(\omega) \leq b_{ij}^k(\omega) y_{ij}, & \forall (i,j) \in \mathcal{A}, \forall k \in \mathcal{K}, \forall \omega \in \Omega, \\
\label{eq:glob_const3_sto} & x_{ij}^k(\omega) \geq 0, & \forall (i,j) \in \mathcal{A}, \forall k \in \mathcal{K}, \forall \omega \in \Omega, \\
\label{eq:glob_const4_sto} & y_{ij} \in \lbrace 0,1 \rbrace, & \forall (i,j) \in \mathcal{A},
\end{alignat} 
where
\[
w_i^k(\omega) =
\begin{cases}
d^k(\omega), & \text{if } i = O(k), \\
-d^k(\omega), & \text{if } i = D(k), \\
0, & \text{otherwise.}
\end{cases}
\]

\subsection{The generator} \label{sec:stoch_gen_tool}

\rev{The description follows the structure of Section~\ref{sec:det_gen_tool}: Section~\ref{sec:stoch_design_workflow} outlines the design and workflow of the stochastic generator, and Section~\ref{sec:stoch_use} describes its use, the accompanying library and its performance.}

\subsubsection{Design and workflow.} \label{sec:stoch_design_workflow}

\rev{The generator reads a base deterministic instance in the STD or DOW format of the deterministic generator (other formats can be converted with simple scripts; see Section~\ref{sec:stoch_inst_real_topo}) and the subsets of parameters to randomize among (i) volumes, (ii) overall arc capacities, (iii) commodity-specific arc capacities and (iv) variable costs. The matrices of target moments and target correlations can either be read in extenso from files, or be generated internally, based on simplifying assumptions, from a few parameters.}

\rev{Target moments can be generated according to one of two distributional characterizations. If the }\emph{uniform} distribution is specified, then each individual randomized parameter is assumed to have first four moments equal to those of a $\textit{uniform}(a, b)$ distribution with $a = \mathcal{D} - (\alpha \cdot \mathcal{D})$, $b = \mathcal{D} + (\beta \cdot \mathcal{D})$, \rev{where $\mathcal{D}$ is the value taken by the parameter in the base deterministic MCFNDP instance and $\alpha$, $\beta$ are spread parameters set by the user satisfying the conditions $\alpha \in [0, 0.95]$, $\beta \in [0, 1000.0)$.} If the \emph{triangular} distribution is specified instead, then each randomized individual parameter is assumed to have first four moments equal to those of a $\textit{triangular}(a, b, c)$ distribution with  $a = \mathcal{D} - (\alpha \cdot \mathcal{D})$, $b = \mathcal{D} + (\beta \cdot \mathcal{D}), c = \mathcal{D}$, \rev{where $\mathcal{D}$, $\alpha$, $\beta$ still have the preceding definitions.} Both distributions are symmetric around $\mathcal{D}$ when $\alpha = \beta$.

\rev{Generated target correlations are specified, in the configuration file or on the command line, by one common value per block of the correlation matrix. Formally, let $\xi=(\xi_1,\dots,\xi_p)$ denote the vector of all randomized parameters, partitioned into the sets of randomized parameters selected by the user among $D$ (volumes), $A$ (overall arc capacities), $C$ (commodity-specific capacities) and $V$ (variable costs). Let $R$ denote the target correlation matrix where $R_{qr}$ is the correlation between $\xi_q$ and $\xi_r$. $R$ is therefore specified blockwise by
\[
R_{qr}=\begin{cases}1 & \text{if } q=r,\\ \rho_{UW} & \text{if } q\neq r,\ \xi_q\in U,\ \xi_r\in W,\end{cases}
\qquad U,W\in\{D,A,C,V\},\ \rho_{UW}=\rho_{WU},
\]
where the block correlations $\rho_{UW}$ are set by the user (zero by default). For example, $\rho_{DD}$ is the common correlation between the volumes of two distinct commodities, and $\rho_{DA}$ the common correlation between any volume and any overall arc capacity. The resulting matrix must be positive definite. Otherwise its Cholesky decomposition fails and the generator stops with an explanatory message.}

Once all target moments and correlations \rev{have been read or generated}, the algorithm 
described in Section~2.5 of \citet{HoylEtAl2003} is supplied with the matrices of target moments and correlations.
From a starting set of scenarios (determined \rev{pseudo-randomly} or fed back from a previous run), the algorithm iteratively applies cubic transformations (to alter moments) and Cholesky factor transformations (to alter correlations) until \rev{the root-mean-square error (RMSE) between realized and target values of the standardized moments and the RMSE between realized and target correlations are both below tolerances set by the user (options \texttt{-m} and \texttt{-c}, both equal to 0.001 by default). Moments are computed with sample estimators. If the algorithm does not meet these tolerances within the prescribed number of iterations, a message is displayed and a new trial starts from new starting scenarios, up to a maximum number of trials set by the user. If the last trial also fails, the generator then displays a diagnosis with remedies, which are discussed in the user guide: checking that the target moments are compatible with a probability distribution and that the target correlation matrix is positive definite, or increasing the numbers of trials, iterations or scenarios. In particular, the condition that the number of scenarios should exceed the number of randomized elements is usually necessary, though not sufficient, for convergence.}

Feasibility of each scenario returned by the HKW algorithm is tested by verifying with either the CPLEX commercial solver or the HiGHS free, open-access solver \citep{HIGHSSite, HuanHall2018} if a solution to the second-stage problem exists under the given scenario when all arcs are open. \rev{The solver is selected when the generator is compiled. Both perform equally well for this check, so that HiGHS is sufficient. The generator also verifies that scenarios are valid, i.e., that their randomized elements are nonnegative, and reports the feasibility and validity of all scenarios. Infeasible scenarios can result from a combination of tightness between overall capacity and commodity volume in the base deterministic instance and high values for the spread parameters $\alpha$ and $\beta$. Invalid scenarios can result from a large value of the spread parameter $\alpha$. Instances comprising infeasible and/or invalid scenarios cannot be repaired by rejecting the defective scenarios since the rejections open a gap between realized and target moments and correlations. A remedy is to choose a base instance with a lesser tightness between volumes and capacities and to choose smaller values for one or both of the spread parameters $\alpha$ and $\beta$.} The generator does not ensure relatively complete recourse of the generated problem instances (i.e., that for every set of admissible first stage values, all scenarios be feasible).

\rev{The generated stochastic instance is} written to the specified output file. In this file, complete deterministic-like problem instances, one for each scenario, are superposed and written under the same format as that used to supply the base deterministic input instance. \rev{Each individual instance associated with a scenario is marked with annotations reporting on its feasibility and validity and stating its probability. This format repeats the non-randomized parameters for every scenario, but it lets the user read each scenario as a complete deterministic instance.}

\subsubsection{Use and performance.} \label{sec:stoch_use}

\rev{The generator, its user guide (USERGUIDE\_STOCH.md) and a commented example (TUTOREXAMPLE\_STOCH.md) are available at \url{https://github.com/larseeri/stoch-MCFNDP-gener}. Like the deterministic generator, it can be built on Linux or run through a Docker image. Figure~\ref{fig:two_commands} shows how the two generators are used in sequence to produce a stochastic instance.}

\rev{
\begin{figure}[htbp]
\centering
\fbox{\begin{minipage}{0.94\textwidth}
\fontsize{9}{11.5}\selectfont\color{blue}
\textbf{Step 1: base deterministic instance} (folder \texttt{program/} of the deterministic generator)\\[2pt]
\texttt{./determ\_gen +F determ\_config\_oneSrcOneSnk\_12oct2026-17h09.par}\\[2pt]
This writes the base instance \texttt{determ\_oneSrcOneSnk\_12oct2026-17h09.std}: a flat 3 x 4 grid (12 nodes) with 40 arcs and 3 commodities, each with one source and one sink. The same instance is also written in the DOW, LP, MPS and TeX formats.
\medskip \\
\textbf{Step 2: stochastic instance} (folder \texttt{program/} of the stochastic generator, where the same base instance is supplied)\\[2pt]
\texttt{./stoch\_gen +F stoch\_config\_oneSrcOneSnk\_12oct2026-17h09\_rstd.par}\\[2pt]
Main options of the configuration file :\\
(at most one option per line in file but comma-separated here; \# precedes comments)\\[2pt]
\begin{tabular}{@{}p{0.43\linewidth}p{0.54\linewidth}@{}}
\texttt{F R}, \texttt{I aboveBaseInstance} & \# base instance of Step 1, read in RSTD format (one source and one sink per commodity)\\
\texttt{S 3} & \# randomize volumes (1) and overall arc capacities (2): 43 random elements\\
\texttt{G}, \texttt{T U}, \texttt{A 0.35}, \texttt{B 0.35} & \# target moments generated from a uniform distribution on $[0.65D, 1.35D]$ around each base value $\mathcal{D}$\\
\texttt{K}, \texttt{XDD 0.5}, \texttt{XAA 0.5}, \texttt{XDA -0.2} & \# target correlations between blocks\\
\texttt{n 500}, \texttt{r 12345}, \texttt{R 200} & \# number of scenarios, seed and stream of the pseudo-random generator\\
\end{tabular}
\medskip \\
The generator writes the target moments and correlations, the 500 raw scenarios, and the stochastic instance \texttt{stoch\_oneSrcOneSnk\_12oct2026-17h09\_rstd.txt}. This file stacks 500 complete deterministic instances, one per scenario, each preceded by a line stating its number, feasibility, validity and probability. The two commands take a few seconds on a laptop, and the files they produce are identical to those of the folder \texttt{examples/oneSrcOneSnk/RSTD} of the stochastic generator's repository.
\end{minipage}}
\caption{\rev{Generating a stochastic instance in two commands}}
\label{fig:two_commands}
\end{figure}
}

\rev{We created with the generator a demonstration library of 416 configurations of stochastic MCFNDP instances, stored in the /library folder of the online repository. The document DESCLIBRARY\_STOCH.txt available at this location details its content. The library illustrates how a moderately diversified set of stochastic instances can be built from a few base instances for a given set of experiments. A stochastic library cannot cover all the settings of the generator for a large set of base instances: the number and size of the base instances, the number of settings and the cost of some of them (many scenarios, many randomized elements) would require prohibitive computing time and storage.}

\rev{The generator is fast enough for its intended low-repetitions use: the 416 instances of the library are built as follows on small to moderately large base instances: Volumes, overall arc capacities or both are randomized in 128 configurations each. Volumes and commodity-specific capacities are randomized in the other 32, built on the smallest base instances. Variable costs are not randomized. Between 1,500 and 6,000 scenarios are generated, with tolerances of 0.001 for the smaller base instances and 0.01 for the larger ones. The largest and median generation times are less than 13 minutes and less than 3 seconds, respectively. Verifying the feasibility of each scenario, a task distinct from generation, is the most time-consuming operation. Generation itself slows down when commodity-specific capacities or variable costs are randomized on networks with many arcs and commodities: the HKW implementation uses non-optimized matrix routines, whose costs grow quickly with the size of the correlation and moment matrices. We return to this limitation in the Conclusion.}

\rev{Table~\ref{tab:scen_comput_verif} reports, for four instances of the historical R series \citep{Canad}, from r05.3 (10 nodes, 60 arcs, 25 commodities) to r18.3 (20 nodes, 315 arcs, 200 commodities), the generation and feasibility-verification times for several numbers of scenarios and two correlation structures, together with the matching errors achieved. (Computations are performed with the default settings of CPLEX 22.1.1 on a single core of an AMD EPYC 7763 processor.) Commodity volumes and overall arc capacities are randomized, with the target moments of a uniform distribution on $[0.75D, 1.25D]$ around each base value $\mathcal{D}$ ($\alpha = \beta = 0.25$). Tolerances about the RMSE made with respect to standardized moments and correlations were both equal to their default value, 0.001, and the numbers of trials and iterations were left at their default values (10 and 50), except for instance r05.3 with 100 scenarios and nonzero correlations, which needed 200 iterations per trial (option \texttt{-i}). Column \emph{RMS corr.} reports the RMSE achieved on correlations, as displayed by the generator. It ranges from 8.18e-5 to 9.96e-4. The RMSE achieved on standardized moments is zero in every run: each iteration of the HKW algorithm ends with the cubic transformation, which matches the four target moments of each randomized parameter exactly, so that convergence is determined by the tolerance on correlations. The notes of the table define its columns. We also produced the content of Table~1 under the alternative assumption of triangular distributions. We found out that average times are similar and that final RMSEs have similar magnitudes.}

\begin{table}[htbp]
\centering
\small
\resizebox{\columnwidth}{!}{

\begin{tabular}{|l|r|r|r|r|r|r|r|r|r|}
\hline
instance & rand. elems & requ. scen. & corr. DD & corr. AA & corr. AD & gener. (sec) & check (sec) & feas. scen. & \rev{RMS corr.} \\
\hline
\hline
r05.3    & 85          & 10000     & 0.0     & 0.0     & 0.0     & \rev{2.0} & \rev{45.7} & \rev{10000} & \rev{\num{8.18e-5}} \\
r05.3    & 85          & 5000      & 0.0     & 0.0     & 0.0     & \rev{1.1} & \rev{22.1} & \rev{5000} & \rev{\num{2.37e-4}} \\
r05.3    & 85          & 1000      & 0.0     & 0.0     & 0.0     & \rev{0.3} & \rev{4.5} & \rev{1000} & \rev{\num{2.67e-4}} \\
r05.3    & 85          & 500       & 0.0     & 0.0     & 0.0     & \rev{0.2} & \rev{2.4} & \rev{500} & \rev{\num{8.71e-4}} \\
r05.3    & 85          & 200       & 0.0     & 0.0     & 0.0     & \rev{0.2} & \rev{1.1} & \rev{200} & \rev{\num{8.69e-4}} \\
r05.3    & 85          & 100       & 0.0     & 0.0     & 0.0     & \rev{0.6} & \rev{0.5} & \rev{100} & \rev{\num{9.69e-4}} \\
\hline
r05.3    & 85          & 10000     & 0.5     & 0.5     & -0.2    & \rev{4.0} & \rev{46.6} & \rev{10000} & \rev{\num{1.09e-4}} \\
r05.3    & 85          & 5000      & 0.5     & 0.5     & -0.2    & \rev{2.0} & \rev{22.0} & \rev{5000} & \rev{\num{1.09e-4}} \\
r05.3    & 85          & 1000      & 0.5     & 0.5     & -0.2    & \rev{0.4} & \rev{4.4} & \rev{1000} & \rev{\num{3.33e-4}} \\
r05.3    & 85          & 500       & 0.5     & 0.5     & -0.2    & \rev{0.3} & \rev{2.3} & \rev{500} & \rev{\num{7.51e-4}} \\
r05.3    & 85          & 200       & 0.5     & 0.5     & -0.2    & \rev{0.2} & \rev{0.9} & \rev{200} & \rev{\num{7.46e-4}} \\
r05.3    & 85          & 100       & 0.5     & 0.5     & -0.2    & \rev{18.9} & \rev{0.5} & \rev{100} & \rev{\num{9.96e-4}} \\
\hline
\hline
r09.3    & 133         & 10000     & 0.0     & 0.0     & 0.0     & \rev{3.8} & \rev{51.3} & \rev{10000} & \rev{\num{1.18e-4}} \\
r09.3    & 133         & 5000      & 0.0     & 0.0     & 0.0     & \rev{2.0} & \rev{27.6} & \rev{5000} & \rev{\num{3.50e-4}} \\
r09.3    & 133         & 1000      & 0.0     & 0.0     & 0.0     & \rev{0.8} & \rev{5.5} & \rev{1000} & \rev{\num{4.82e-4}} \\
r09.3    & 133         & 500       & 0.0     & 0.0     & 0.0     & \rev{0.6} & \rev{2.7} & \rev{500} & \rev{\num{4.82e-4}} \\
r09.3    & 133         & 200       & 0.0     & 0.0     & 0.0     & \rev{0.6} & \rev{1.3} & \rev{200} & \rev{\num{8.59e-4}} \\
\hline
r09.3    & 133         & 10000     & 0.5     & 0.5     & -0.2    & \rev{7.7} & \rev{55.0} & \rev{10000} & \rev{\num{1.07e-4}} \\
r09.3    & 133         & 5000      & 0.5     & 0.5     & -0.2    & \rev{3.9} & \rev{25.5} & \rev{5000} & \rev{\num{1.16e-4}} \\
r09.3    & 133         & 1000      & 0.5     & 0.5     & -0.2    & \rev{0.9} & \rev{5.4} & \rev{1000} & \rev{\num{4.58e-4}} \\
r09.3    & 133         & 500       & 0.5     & 0.5     & -0.2    & \rev{0.6} & \rev{2.7} & \rev{500} & \rev{\num{4.63e-4}} \\
r09.3    & 133         & 200       & 0.5     & 0.5     & -0.2    & \rev{0.8} & \rev{1.1} & \rev{200} & \rev{\num{9.01e-4}} \\
\hline
\hline
r14.3    & 320         & 10000     & 0.0     & 0.0     & 0.0     & \rev{14.9} & \rev{117.9} & \rev{10000} & \rev{\num{2.59e-4}} \\
r14.3    & 320         & 5000      & 0.0     & 0.0     & 0.0     & \rev{7.6} & \rev{57.5} & \rev{5000} & \rev{\num{6.93e-4}} \\
r14.3    & 320         & 1000      & 0.0     & 0.0     & 0.0     & \rev{4.4} & \rev{12.0} & \rev{1000} & \rev{\num{4.97e-4}} \\
r14.3    & 320         & 500       & 0.0     & 0.0     & 0.0     & \rev{5.1} & \rev{6.0} & \rev{500} & \rev{\num{7.97e-4}} \\
\hline
r14.3    & 320         & 10000     & 0.5     & 0.5     & -0.2    & \rev{33.3} & \rev{115.9} & \rev{10000} & \rev{\num{1.40e-4}} \\
r14.3    & 320         & 5000      & 0.5     & 0.5     & -0.2    & \rev{17.1} & \rev{58.9} & \rev{5000} & \rev{\num{2.32e-4}} \\
r14.3    & 320         & 1000      & 0.5     & 0.5     & -0.2    & \rev{5.6} & \rev{11.8} & \rev{1000} & \rev{\num{5.86e-4}} \\
r14.3    & 320         & 500       & 0.5     & 0.5     & -0.2    & \rev{6.6} & \rev{6.0} & \rev{500} & \rev{\num{7.78e-4}} \\
\hline
\hline
r18.3    & 515         & 10000     & 0.0     & 0.0     & 0.0     & \rev{35.3} & \rev{245.8} & \rev{10000} & \rev{\num{3.88e-4}} \\
r18.3    & 515         & 5000      & 0.0     & 0.0     & 0.0     & \rev{18.3} & \rev{124.3} & \rev{5000} & \rev{\num{9.60e-4}} \\
r18.3    & 515         & 1000      & 0.0     & 0.0     & 0.0     & \rev{15.8} & \rev{25.1} & \rev{1000} & \rev{\num{6.40e-4}} \\
\hline
r18.3    & 515         & 10000     & 0.5     & 0.5     & -0.2    & \rev{86.1} & \rev{232.0} & \rev{10000} & \rev{\num{1.95e-4}} \\
r18.3    & 515         & 5000      & 0.5     & 0.5     & -0.2    & \rev{42.4} & \rev{122.7} & \rev{5000} & \rev{\num{3.95e-4}} \\
r18.3    & 515         & 1000      & 0.5     & 0.5     & -0.2    & \rev{19.1} & \rev{24.4} & \rev{1000} & \rev{\num{9.66e-4}} \\
\hline
\multicolumn{10}{l}{\small \rev{instance: base deterministic instance of the R series.}}\\
\multicolumn{10}{l}{\small \rev{rand. elems: number of parameters randomized between scenarios.}}\\
\multicolumn{10}{l}{\small \rev{requ. scen.: number of scenarios requested.}}\\
\multicolumn{10}{l}{\small \rev{corr. DD: common target cross-correlation coefficient between individual commodity volumes.}}\\
\multicolumn{10}{l}{\small \rev{corr. AA: common target cross-correlation coefficient between individual arc capacities.}}\\
\multicolumn{10}{l}{\small \rev{corr. AD: common target cross-correlation coefficient between individual arc capacities and commodity volumes.}}\\
\multicolumn{10}{l}{\small \rev{gener. (sec): scenario generation time, in seconds.}}\\
\multicolumn{10}{l}{\small \rev{check (sec): feasibility verification time, in seconds.}}\\
\multicolumn{10}{l}{\small \rev{feas. scen.: number of feasible scenarios.}}\\
\multicolumn{10}{l}{\small \rev{RMS corr.: RMSE between realized and target correlations.}}\\
\end{tabular}
}
\caption{Scenario computation and verification times\rev{, and moment-matching errors}}
\label{tab:scen_comput_verif}
\end{table}

\subsection{\rev{Stochastic instances on benchmark hub-and-spoke networks}}
\label{sec:stoch_inst_real_topo}

\rev{The stochastic generator is not limited to base instances produced by our deterministic generator. It accepts any deterministic MCFNDP stored in the DOW or STD format. To illustrate this, we built a set of 20 two-stage stochastic MCFNDPs on the hub-and-spoke topologies of \citet{BonnEtAl2025}, which aim to reproduce characteristics of real networks and are designed to match targets for complex-network metrics rather than being drawn from elementary ranges. The base instances, the conversion script, the configuration of every generation run and a detailed description of the set are available at \url{https://github.com/larseeri/stoch-MCFNDP-gener/tree/main/external\_sources/Bonnet-et-al\_2025}.}

\rev{\citet{BonnEtAl2025} provide 180 benchmark networks with 50 nodes, obtained by combining four arc densities (5, 20, 50 and 80\%, that is, 122 to 1960 arcs), three hub node ratios (1, 10 and 20\%) and three levels of reciprocity (20, 50 and 90\%), with five replicates of each combination. The level of reciprocity indicates the degree to which arcs between pairs of nodes exist in both directions. Every arc has a capacity of 1,000. We selected 20 networks to obtain a good balance in terms of arc density, hub node ratio and level of reciprocity. As the benchmark networks contain no commodities, the volumes were generated with the code of \citet{BonnEtAl2025} itself, with its default settings for static instances: origin--destination pairs are drawn among the pairs joined by a path, and volumes follow a normal distribution with mean 100 and standard deviation 50, truncated to $[10, 1000]$. Each network received 50, 100 or 200 commodities, with its own seed. Our repository gives the design and a script that regenerates the base instances with this code. Appendix B also provides more details on the scenario generation process.}

\section{Conclusion} \label{sec:concl}

\rev{We introduced two fast and flexible generators of deterministic and stochastic MCFNDP instances, which we hope will facilitate systematic experimentation and foster the reproducibility and comparability of published research. The deterministic generator modernizes, extends and documents Mulgen. It} stands on its own and also supplies base instances to the stochastic generator. 

\rev{The stochastic generator, which is new, randomizes a user-specified subset of the parameters of any deterministic instance supplied in an admissible format, using a moment-matching algorithm driven by a few interpretable requirements. Both come with libraries of configurations from which 6,000 deterministic and 416 stochastic instances can be regenerated.}

As a natural extension of this line of research, we envision the creation of similar generators simulating instances of other standardized problem families of high interest, \rev{such as facility location, hub location, and location-routing problems}. \rev{Evidence collected across several problem families would make the assessment of solution methods more reliable.}

\rev{Two limitations also suggest directions for future work. \emph{First}, the implementation of the HKW algorithm relies on non-optimized matrix routines, so that generation slows down markedly when commodity-specific capacities or variable costs are randomized on large networks. An optimized linear-algebra library such as OpenBLAS \citep{openblas} would alleviate this bottleneck. \emph{Second}, the deterministic generator draws its parameters from uniform distributions over user-specified ranges, which reflects diffuse prior knowledge. Other distributions, as well as other target distributions in the stochastic generator, could easily be added.}

\section*{Acknowledgments}
This research was funded by the Canadian National Railway Company (CN) Chair in Optimization of Railway Operations at Université de Montréal, the Canada Research Chair program [950-232244] and IVADO fundamental research program on integrated machine learning and optimization for decision making under uncertainty. We thank Antonio Frangioni for his stewardship of the Canad instances and the Mulgen generator on the CommaLAB website. We are also grateful for important insights provided by verifications and testing conducted by Ammar Metnani. The original code of the Mulgen generator was written by Bernard Gendron who passed away at a too young age. We are deeply grateful for the invaluable discussions on network design we had the privilege to have with him. Finally, we thank the Associate Editor and three anonymous reviewers for their valuable comments which helped improve the quality of this paper.


\section*{Appendix \rev{A}. Existing methods for the construction of scenarios} 

We \rev{provide an overview of} the methods proposed in the literature on linear two-stage stochastic programming for constructing scenarios. We aim to identify those that are suited for the purpose of synthesizing a distribution and, among the latter, discuss usability in our particular context of application. Detailed reviews of the literature on scenario construction for the purposes of stochastic programming are available in \citet{Lohn2016, BounEtAl2022, Kaut2021, KeutEtAl2023, RujeEtAl2022, BertMund2022}. We distinguish the following operating principles and classify the existing methods accordingly.

\paragraph{Sampling scenarios.}
Surveys of the literature on random and quasi-random sampling of scenarios are available in \citet{Shap2003, Glas2004, BayrEtAl2011, HomeBayr2014}. Later advances are available in \citet{LeovRomi2015, Lohn2016} (with stratification through Voronoi sampling). Methods in this class sample among existing scenarios included in an assumed known probability distribution. They are therefore inadequate for the purpose of synthesizing a distribution based on primitive requirements.

\paragraph{Clustering scenarios.}
Existing scenarios included in an assumed known probability distribution are clustered with k-means, k-medians, k-medoids, either based on values of parameters appearing in scenarios or on values of second stage outcomes associated to individual scenarios \citep{KeutEtAl2023, HewiEtAl2022, AbouEtAl2022}. These methods are inadequate for our purposes as they operate on scenarios of an assumed known probability distribution.

\paragraph{Minimization of probabilistic distance.}
This area of research is  substantial and currently the most active. It also extends to multistage stochastic programming and scenario trees. Those topics exceed our scope. Surveys are available in \citet{Lohn2016, RujeEtAl2022, BertMund2022}. Theoretical advances and methodological overviews are available from \citet{DupaEtAl2003, HeitRomi2003, HeitRomi2007, HenrEtAl2008, HenrEtAl2009, PfluPich2015, HenrRomi2022, RujeEtAl2022}. \rev{Large-scale}, industrial applications are presented in \citet{LiFlou2014, LiFlou2016, RujeEtAl2022, BertMund2022, AbouEtAl2022}. Methods in this class aim to minimize or show convergence of the Wasserstein (a.k.a. Kantorovich-Rubinstein) or Fortet-Mourier probabilistic distances or an upper bound thereof between a known, exact distribution and a new, approximate distribution resulting from a constructed set of scenarios. There is a partial overlap between this class and the preceding ones as sampling and clustering may be involved in building the sets of scenarios aimed at controlling the probabilistic distance \citep[see for instance][]{DupaEtAl2003, HeitRomi2003, HeitRomi2007, HenrEtAl2008, HenrEtAl2009, HenrRomi2022, RujeEtAl2022, BertMund2022, AbouEtAl2022}.
Methods in this class control the probabilistic distance either (i) between the exact and approximate distributions of the second stage parameters while relying on stability analysis to characterize the probabilistic distance between the resulting outputs, or (ii) directly between the distributions of the outputs achieved under the exact and approximate distributions of the parameters. Whereas it would be conceivable to devise constructive methods based on primitive requirements aiming to minimize probabilistic distances, this would exceed the scope of our proposition.
 
\paragraph{Matching moments, correlations and distribution functions.}
Methods in this class build scenarios by matching target moments and correlations or matching target distributions for the stochastic parameters \citep[see][]{Flei1978, HoylWall2001, HoylEtAl2003, MehrPapp2013}.
This type of method is applicable for our purposes as it is conceivable to state target moments and correlations for the stochastic parameters of the MCFNDP. \rev{An efficient implementation by one of the authors of \citet{HoylEtAl2003} has been used to build scenarios for MCFNDPs in several publications \citep[see, e.g.,][]{crainic2011progressive,CraiEtAl2021b}.} \rev{This implementation} requires supplying targets about first four moments and linear correlations for all stochastic parameters in the MCFNDP at hand. We saw in Section~\ref{sec:stoch_design_workflow} how this information may be derived from a small set of assumptions. Using the more general method proposed in \citet{MehrPapp2013} is also conceivable. In contrast with \citet{HoylEtAl2003}, \citet{MehrPapp2013} can match the distribution function rather than only first four moments and linear correlations. However, its greater informational prerequisites are difficult to postulate convincingly in practice and its computational tractability stands at a prototypical stage. These reasons lead us to prefer the method of \citet{HoylEtAl2003} to that of \citet{MehrPapp2013}. \citet{Kaut2021} proposes a MILP solution to the moment-matching problem that might be viewed as an alternative to the iterative NLP algorithm of \citet{HoylEtAl2003}. However, we prefer the latter as the former is prototypical. Moment-matching methods can also be joined with clustering methods \citep{LiZhu2016}, principal component analysis \citep{ChopSelv2020}, distribution matching \citep{CalfEtAl2014, BounEtAl2022}. These variants require a detailed knowledge of underlying scenarios and distribution that is unavailable in our context of application.

\paragraph{Copula sampling.}
A number of scenario selection methods account for dependencies in the assumed available underlying joint distribution through copula sampling and may combine it with matching of moments, matching of marginal distributions or clustering to account for marginal distributions \citep[see][]{SutiPran2007, KautWall2011, Kaut2014, Kaut2015, QiuEtAl2019}.
In our context of application, while it would be conceivable to select a copula based on primitive assumptions and conduct sampling from it, this would presuppose a level of knowledge about the dependencies among the stochastic parameters of the MCFNDP that may be difficult to justify in practice.

\section*{\rev{Appendix B. Details of the scenario generation process for the Bonnet et al.~instances}}

\rev{The instance files of \citet{BonnEtAl2025} list nodes, arcs (unit cost, fixed cost and capacity) and commodities (origin, destination and volume) in a format of their own. A Bash script, \texttt{bonnet2std.sh}, available in our repository, converts them to the DOW or the STD format read by the stochastic generator. It renumbers nodes from 1 and keeps the order of arcs and commodities. Since these instances have no commodity-specific data, the STD output assigns to every commodity the unit cost and the capacity of the arc. Our set uses the DOW format, in which values are copied unchanged. The script can also write the STD format, with values rounded to integers by default or copied unchanged with option \texttt{-x}.}

\rev{Each base instance was then randomized with the stochastic generator. Only the volumes are random (\texttt{-S 1}). Costs, capacities and origin--destination pairs are those of the base instance. Target moments are those of a uniform distribution on $[(1-s)\mathcal{D}, (1+s)\mathcal{D}]$ around each base volume~$\mathcal{D}$ ($\alpha=\beta=s$), and distinct volumes have a common target correlation~$\rho$ (\texttt{-K -XDD $\rho$}). Each of the ten combinations of $s \in \{0.1, 0.2, \ldots, 0.5\}$ and $\rho \in \{0, 0.5\}$ was assigned to two base instances. Each instance with $K$ commodities has $2K$ equiprobable scenarios (100 to 400), and the HKW algorithm \citep{HoylEtAl2003} was run with a tolerance of 0.001 on moments and correlations. We required all scenarios of an instance to be feasible and valid. Whenever the generator reported an infeasible scenario, the attempt was rejected and $s$ was lowered by 0.05 for a new attempt. Every run of the generator is a numbered experiment~$e$, which uses the seed 50163844 (\texttt{-r}) and the stream $47 + 100(e-1)$ (\texttt{-R}). Nineteen base instances were accepted at their first attempt. The twentieth, whose capacities are the tightest of the set (its volumes can be scaled by at most 1.40 and still be routed with all arcs open), had 6 and 1 infeasible scenarios out of 200 with $s = 0.5$ and $s = 0.45$, and none with $s = 0.4$, for a total of 22 experiments.}

\rev{In all accepted instances, the HKW algorithm met its tolerance in its first trial. The mean ratio of scenario volume to base volume is 1.000, its standard deviation is within 0.6\% of its target $s/\sqrt{3}$, and the average correlation between distinct volumes lies between 0.4994 and 0.4998 for $\rho = 0.5$ and is 0.000 for $\rho = 0$. As the scenario files total 620~MB, we do not store them in the repository. Instead, we store the configuration file read by the generator, which holds the content of its command line with the seed and the stream, together with its log, its outcome and its generation time, as well as the target moments and correlations and the realized statistics of each instance (about 4~MB in total).}

\bibliographystyle{plainnat_custom}

\bibliography{refs}

\end{document}